\documentclass[11pt]{amsart}
\usepackage{amssymb}

\newtheorem*{thm}{Theorem}

\newcommand\F{\mathbb{F}}

\begin{document}


\title{On a theorem of Carlitz}

\author{Michael E. Zieve}
\address{
  Department of Mathematics,
  University of Michigan,
  Ann Arbor, MI 48109--1043,
  USA
}
\email{zieve@umich.edu}
\urladdr{www.math.lsa.umich.edu/$\sim$zieve/}


\begin{abstract}
Carlitz proved that, for $q>2$, the group of all permutations of $\F_q$
is generated by the permutations induced by degree-one polynomials and
$x^{q-2}$.  His proof relies on a remarkable polynomial which appears
to have been found by magic.  We show here that no magic is required:
there is a straightforward
way to produce a simple polynomial which has the same remarkable properties
as the complicated polynomial in Carlitz's proof.  We also identify the
crucial subtlety which allows such simple polynomials to exist, and discuss
some consequences.
\end{abstract}

\maketitle


The theorem in the title is as follows:

\begin{thm}
If $q>2$ then every permutation of\/ $\F_q$ is the composition of permutations
induced by $x^{q-2}$ and by degree-one polynomials over\/ $\F_q$.
\end{thm}

Betti proved this for $q=5$ (as the final assertion in \cite{Betti}),
and Dickson proved it for $q=7$ \cite[p.~119]{Dicksonthesis}.
In response to a question posed by Straus, Carlitz proved
it in general \cite{Carlitz}, via the following argument.  It suffices to prove the
result in case the permutation is a $2$-cycle of the form $(0a)$ with
$a\in\F_q^*$, since
every permutation is a product of such $2$-cycles.  And Carlitz observed
that $(0a)$ is the permutation of $\F_q$ induced by
\[
f_a(x) := -a^2\Bigl(\Bigl((x-a)^{q-2}+\frac{1}{a}\Bigr)^{q-2}-a\Bigr)^{q-2}.
\]
Although it is straightforward to verify that $f_a$ induces the permutation $(0a)$,
it is not at all clear how one could have discovered the polynomial
$f_a$ in the first place.  Indeed, several authors have presented $f_a$
as a mysterious and complicated object: for instance, \cite[p.~169]{Wells}
and \cite[p.~358]{LN} assert that
this representation of $(0a)$ demonstrates that ``simplicity as
polynomials and simplicity as permutations are not equivalent."

My purpose here is to remove the mystery from Carlitz's proof,
by presenting a straightforward procedure for producing a
simple polynomial which has the same crucial property as $f_a$, namely
that of inducing the permutation $(0a)$.
Note that $\mu(x):=1-1/x$
induces an order-$3$ permutation of $\F_q\cup\{\infty\}$, and one cycle
of $\mu$ is $(\infty 1 0)$.  Then $h(x):=1-x^{q-2}$ agrees with $\mu$ on
$\F_q^*$, and $h$ interchanges $0$ and $1$, so
$g(x):=h(h(h(x)))$ induces the permutation $(01)$ on $\F_q$.
Thus $a g(x/a)$ induces the permutation $(0a)$.

The surprising feature of this proof -- and of Carlitz's result, once we
identify $x^{q-2}$ with $1/x$ -- is that we have expressed each element of
the symmetric group $S_q$ as a composition
of degree-one rational functions, which should not be possible since the set $G$ of
degree-one rational functions is closed under composition and $\#G=q^3-q$ is
typically much
smaller than $q!$.  However, the two notions of composition are incompatible in a
subtle way, since we are not viewing an action of $G$:
although we begin with the action of $G$ on $\F_q\cup\{\infty\}$,
we must identify $\infty$ with $0$ in order to view $1/x$ as a permutation of
$\F_q$, but we cannot make a compatible identification for degree-one polynomials.

Note that $x^{q-2}$ permutes $\F_{q^k}$ for infinitely many $k$:
specifically, for all $k$ such that $q^k-1$ is coprime to $q-2$, which amounts
to requiring that 
$k$ is not divisible by any of the numbers $r_{\ell}$,
where $\ell$ is a prime factor of $q-2$ and $r_{\ell}$ is the order of $q$
in $\F_{\ell}^*$.  Thus, as noted by Carlitz
\cite[\S 2]{Carlitz2}, any composition of degree-one polynomials and copies of
$x^{q-2}$ will also permute these infinitely many extensions of $\F_q$, so
any such composition is an \emph{exceptional polynomial} (cf.\ \cite{GRZ}
and the references therein).  Hence Carlitz's result implies that every
permutation of $\F_q$ is induced by an exceptional polynomial.

The condition $q>2$ in the above theorem is needed only because $x^{q-2}$
does not permute $\F_q$ when $q=2$.  When $q=2$, every permutation of $\F_q$
is represented by a degree-one polynomial.

Further results related to the above theorem are given in
\cite{Carlitz3,FFA,Fryer,Fryer2,Wells}.




\begin{thebibliography}{9}
\newcommand{\au}[1]{{#1},}
\newcommand{\ti}[1]{\textit{#1},}
\newcommand{\jo}[1]{{#1}}
\newcommand{\vo}[1]{\textbf{#1}}
\newcommand{\yr}[1]{(#1),}
\newcommand{\ppx}[1]{#1,}
\newcommand{\pp}[1]{#1.}
\newcommand{\pps}[1]{#1}
\newcommand{\bk}[1]{{#1},}
\newcommand{\inbk}[1]{in: \bk{#1}}
\newcommand{\xxx}[1]{{arXiv:#1}.}


\bibitem{Betti}
\au{E. Betti}
\ti{Sopra la risolubilit\'a per radicali delle equazioni algebriche
irriduttibili di grado primo}
\jo{Annali Sci. Mat. Fis.}
\vo{2}
\yr{1851}
\pp{5--19}
[= Op. Mat.
\vo{1}
\yr{1903}
\pp{17--27}]

\bibitem{Carlitz}
\au{L. Carlitz}
\ti{Permutations in a finite field}
\jo{Proc. Amer. Math. Soc.}
\vo{4}
\yr{1953}
\pp{538}

\bibitem{Carlitz2}
\au{\bysame}
\ti{Permutations in finite fields}
\jo{Acta Sci. Math. (Szeged)}
\vo{24}
\yr{1963}
\pp{196--203}

\bibitem{Carlitz3}
\au{\bysame}
\ti{A note on permutations in an arbitrary field}
\jo{Proc. Amer. Math. Soc.}
\vo{14}
\yr{1963}
\pp{101}

\bibitem{FFA}
\au{A. \c Ce\c smelio\u glu, W. Meidl and A. Topuzo\u glu}
\ti{On the cycle structure of permutation polynomials}
\jo{Finite Fields Appl.}
\vo{14}
\yr{2008}
\pp{593--614}

\bibitem{Dicksonthesis}
\au{L.~E. Dickson}
\ti{The analytic representation of substitutions on a power of a prime
number of letters with a discussion of the linear group}
\jo{Ann. of Math.}
\vo{11}
\yr{1896}
\pp{65--120}

\bibitem{Fryer}
\au{K.~D. Fryer}
\ti{A class of permutation groups of prime degree}
\jo{Canad. J. Math.}
\vo{7}
\yr{1955}
\pp{24--34}

\bibitem{Fryer2}
\au{\bysame}
\ti{Note on permutations in a finite field}
\jo{Proc. Amer. Math. Soc.}
\vo{6}
\yr{1955}
\pp{1--2}

\bibitem{GRZ}
\au{R.~M. Guralnick, J.~E. Rosenberg and M.~E. Zieve}
\ti{A new family of exceptional polynomials in characteristic $2$}
\jo{Annals of Math.}
\vo{172}
\yr{2010}
\pp{1367--1396}

\bibitem{LN}
\au{R.~Lidl and H.~Niederreiter}
\bk{Finite Fields}
Addison--Wesley, Reading, 1983.

\bibitem{Wells}
\au{C. Wells}
\ti{Generators for groups of permutation polynomials over finite fields}
\jo{Acta Sci. Math. Szeged.}
\vo{29}
\yr{1968}
\pp{167--176}

\end{thebibliography}
\end{document}